\documentstyle[10pt,twoside]{article}
\input xypic
\xyoption{all}
\LaTeXdiagrams
\xyoption{2cell}
\UseAllTwocells
\newtheorem{defi}{Definition}
\newtheorem{teo}{Theorem}

\newtheorem{prop}{Proposition}

\def\be{\begin{equation}}
\def\ee{\end{equation}}
\def\lra{\longrightarrow}
\def\Lra{\Longrightarrow}
\def\lrat{\tilde{\longrightarrow}}
\def\Lrat{\tilde{\Longrightarrow}}
\def\fin{Fin_{sk}}
\def\finbar{\overline{\bf Fin}_{sk}}
\def\A{{\bf A}}
\def\B{{\bf B}}
\def\C{{\bf C}}
\def\Cat{{\bf Cat}}
\def\D{{\bf D}}

\def\N{{\bf N}}
\def\O{{\cal O}}
\def\X{{\bf X}}

\def\T{{\bf T}}
\def\To{{\bf T^{op}}}%
\newcommand{\TT}[1]{{\bf T}_{#1}}

\def\Tassoc{{\bf T}_{Assoc}}
\def\Tmon{{\bf T}_{Mon}}
\def\Tsmon{{\bf T}_{sMon}}
\def\Tsym{{\bf T}_{Sym}}
\def\Tbin{{\bf T}_{Bin}}
\def\Tbraid{{\bf T}_{Braid}}
\def\Tsbraid{{\bf T}_{sBraid}}
\def\Tcomm{{\bf T}_{Comm}}
\def\Tscomm{{\bf T}_{sComm}}
\def\Tstar{{\bf T}_{Strfsh}}
\def\Tssym{{\bf T}_{sSym}}
\def\Tbal{{\bf T}_{Bal}}
\def\Tx{{\bf T}_{X}}
\def\Tsx{{\bf T}_{sX}}
\def\Tsbal{{\bf T}_{sBal}}
\def\th{{\bf Theories}}
\def\Th{\widetilde{\bf Theories}}
\def\TTh{\widetilde{\bf 2Theories}}
\def\TCat{\widetilde{\bf 2Cat}}
\def\TCatcat{\widetilde{\bf 2Cat/ Cat}}

\def\KP{\otimes^K}
\def\TA{{\bf 2Alg}}
\newenvironment{examp}{ \stepcounter{examnum} {\bf \noindent Example.\arabic{section}.\arabic{examnum}:}}{$\Box$}

\newcounter{examnum}[section]
\newcounter{remarnum}[section]
\begin{document}
\title{Coherence, Homotopy and 2-Theories}
\author{Noson S. Yanofsky\\ Department of Computer and Information Science\\ Brooklyn College, CUNY\thanks{This research was supported by
PSC-CUNY award \#61791-00-30}}
\date{June 1, 2000}
\maketitle
\begin{abstract}
\noindent 2-Theories are a canonical way of describing categories
with extra structure. 2-theory-morphisms are used when discussing
how one structure can be replaced with another structure. This is
central to categorical coherence theory. We place a Quillen model
category structure on the category of 2-theories and
2-theory-morphisms where the weak equivalences are biequivalences
of 2-theories. A biequivalence of 2-theories (Morita equivalence)
induces and is induced by a biequivalence of 2-categories of
algebras. This model category structure allows one to talk of the
homotopy of 2-theories and discuss the universal properties of
coherence.
\end{abstract}
\section{Introduction}
The history of coherence theory has its roots in homotopy theory.
Saunders Mac Lane's foundational paper \cite{MacLane} on coherence
theory was an abstraction of earlier work by James Stasheff on
H-spaces \cite{Stash} and by D.B.A. Epstein on Steenrod
operations. Coherence theory went on to become an important part
of many diverse areas of computer science and mathematics.
Questions of categorical coherence arise in, to name but a few
areas, linear logic, proof theory, concurrency theory,
low-dimensional topology, quantum groups and quantum field theory.
Although coherence theory has become a mature and independent
part of category theory, it has always had a distinct homotopy
theory flavor. One has the feeling that a monoidal category is
the same ``up to homotopy'' as a strict monoidal category. Or
that a braided tensor category can be ``deformed'' to a strict
braided tensor category. This paper is a step toward clarifying
and formulating the exact relationship between coherence theory
and homotopy theory.

We shall use the language of algebraic 2-theories to talk about
coherence. Algebraic 2-theories are a generalization of Lawvere's
algebraic (1-)theories which are central to his functorial
semantics \cite{Lawvere}.  Algebraic 1-theories are a categorical
description of sets (or topological spaces, manifolds, vector
spaces, etc) with extra equational structure. Algebraic
2-theories are 2-categorical descriptions of categories (or any
objects in a 2-category) with extra structure. Whereas the
1-cells in a 1-theory correspond to operations on sets, the
1-cells in a 2-theory correspond to functors. 2-cells in a
2-theory correspond to natural transformations between functors.
So we have the 2-theory of monoidal categories, $\Tmon$, strict
monoidal categories $\Tsmon$, braided categories $\Tbraid$,
balanced categories $\Tbal$ etc (see \cite{Joy&Str} for many
examples and uses.) Each 2-theory $\T$ has an associated
2-category $\TA(\T,\Cat)$ of algebras, morphisms between algebras
and natural transformations between morphisms.

John Gray was the first to define 2-theories \cite{Grayb}. They
are also discussed in Bloom {\it et al} \cite{Bloom} and John
Powers \cite{Powers1, Powers2} (although our definition of
morphisms between algebras is slightly different from
\cite{Powers2})

Coherence theory is concerned with the way one algebraic structure
interacts with another. Let $F:\T_1 \lra \T_2$ be a
2-theory-morphism. $F$ induces a $F^*:\TA(\T_2,\Cat) \lra
\TA(\T_1,\Cat)$ and a quasi-left adjoint $F_*:\TA(\T_1,\Cat) \lra
\TA(\T_2,\Cat)$. The strength of this adjunction determines, and
is determined by $F$. The adjunction can be an isomorphism,
equivalence, biequivalence, strict quasi-adjunction etc. We shall
be particularly interested in the notion of a biequivalence.

A biequivalence is a 2-categorical generalization of an
equivalence. Whereas in an equivalence, the unit and counit are
isomorphisms, in a biequivalence, the unit and counit are
themselves equivalences. Examples of biequivalences are abundant.
The phrase ``Every monoidal category is tensor equivalent to a
strict monoidal category'' is a way of saying that
$\TA(\Tmon,\Cat)$ and $\TA(\Tsmon, \Cat)$ are biequivalent. We
might call $\Tmon$ and $\Tsmon$ Morita equivalent 2-theories.
This fact is equivalent to saying that $\Tmon$ is biequivalent to
$\Tsmon$. When we say that two types of categorical structure are
the same ``up to homotopy'' or ``are of the same homotopy type'',
we mean there is a biequivalence between their 2-theories.

In 1967, Daniel Quillen \cite{Quillen} showed one how to talk
about homotopy theory from a categorical point of view (see also
\cite{Dwyer&S, Hovey}). A category $\C$ has a functorial closed
Quillen model category (FCQMC) structure if there are three
classes of morphisms in $\C$ called weak equivalences, fibrations
and cofibrations. These classes of maps must satisfy certain
axioms that are of importance in homotopy theory. Examples of
FCQMCs are the category of topological spaces, simplicial sets,
simplicial complexes and $\Cat$. However many other categories
like algebraic objects, chain complexes of algebraic objects, or
spectra also have FCQMC structures. Once one has such a
structure, one can go on to create the analogies of path spaces,
mapping cylinders, Puppe sequences and all other important tools
of homotopy theory. Given a FCQMC structure, one formally inverts
the weak equivalences (that is, makes the weak equivalences into
isomorphisms) in order to construct $Ho(\C)$, the homotopy
category of $\C$, and a functor $\gamma : \C \lra Ho(\C)$.
$\gamma$ has the universal property that given any category $\D$
and functor $F:\C \lra \D$ that inverts the weak equivalences of
$\C$, there is a unique $G:Ho(\C) \lra \D$ such that $G\circ
\gamma = F$.

The goal of this paper is to show that the category of 2-theories
and 2-theory-morphisms has a FCQMC structure. The weak
equivalences will be biequivalences. It is important to realize
that this paper is not a generalization of any known theorem about
1-theories. We know of no nontrivial FCQMC structure on the
category of 1-theories (a FCQMC structure is trivial if the weak
equivalences are exactly isomorphisms and hence $\C=Ho(\C)$.) Only
2-theories have the flexibility to have a homotopy theory.
1-theories are too rigid for this. If one tries to do the same
trick with 1-theories by making weak equivalences into genuine
equivalences of 1-theories, one gets the trivial FCQMC structure.
This follows from the fact that $F$ is a genuine equivalences of
1-theories iff $F$ is an isomorphism of 1-theories.

With the FCQMC in place, we can go on and write down universal
properties of coherence.

This paper is organized as follows. Section 2 is a review of the
relevant aspects of 2-theories. Section 3 is a discussion of the
calculus of biequivalences. After defining the weak equivalences,
fibrations, and cofibrations, Section 4 goes on to prove that
they satisfy the axioms of an FCQMC structure. Some universal
properties of coherence are discussed in Section 5. Section 6 is
a look towards future directions that this project can take.

{\bf Notation.} 2-Theories shall be denoted as $\T, \T', \TT 1,
\TT 2 , \ldots $ 2-theory-morphisms are capital letters $F, G, H,
\ldots$. Lower-case Greek letters $\alpha, \beta,
 \gamma, \delta , \ldots$ will denote 2-theory-natural
 transformations. Capital Greek letters $\Phi, \Psi, \Xi, \ldots$
 will denote 2-theory-modifications. The three compositions of
 morphisms in $\TTh$ shall be denoted $\circ_0, \circ_1,$ and
 $\circ_2$ but will be omitted when no ambiguity arises.
 We shall write $x \in_i \X$ if $x$ is an i-cell of the
 n-category $\X$ where $i=0,1,2, \ldots , n$. 1-Categories
 will be shown in non-bold typeface while 2-categories are in
bold
 typeface. Following Gray \cite{Gray}, we shall denote all
 3-categories by placing a tilde above it.

{\bf Acknowledgments.} I am grateful to Alex Heller for teaching
me all about model categories (and many other subjects) and for
listening to this work while it was in progress.  I am thankful
for helpful conversations with David Blanc, Mirco Mannucci and
Dennis Sullivan.

\section{2-Theories}
Let $\fin$ denote the skeletal category of finite sets. The
2-category $\finbar$ is $\fin$ with only identity 2-cells. Place
a coproduct structure on $\finbar$. A coproduct structure for a
2-category is similar to a coproduct structure for a 1-category.
However, there is an added requirement that for every finite
family of 1-cells with common source and target, there is a
1-cell with injection 2-cells that satisfy the obvious universal
property. When we talk of preserving coproduct structures, we
mean preserving the coproduct strictly (equality).

\begin{defi}
A (single sorted algebraic) {\bf 2-theory} is a 2-category $\T$
 with a given coproduct structure and a 2-functor $G_{\T}:
\finbar \longrightarrow \T$ such that $G_{\T}$ is bijective on
0-cells and preserves the coproduct structure.
\end{defi}

The following examples are well known.

\begin{examp}
$\finbar$ is the initial 2-theory. Just as $\fin$ is
the theory of sets, so too,  $\finbar$ is the theory of
categories.
\end{examp}

\begin{examp}
Let $\Tbin$ be $\finbar$ with a nontrivial generating 1-cell
$\otimes:1 \longrightarrow 2$ thought of as a binary operation
(bifunctor).
\end{examp}

\begin{examp}
$\Tmon$ is the 2-theory of monoidal (tensor) categories. It is a
2-theory ``over'' $\Tbin$ with a 1-cell $e:1 \longrightarrow 0$.
The isomorphic 2-cells are generated by $$
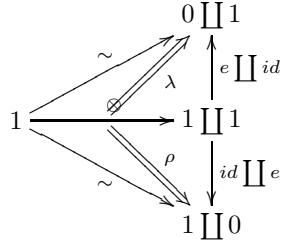
\xymatrix{ && 0 \coprod
1
\\
1 \ar[urr]^{\sim} \ar[rr]^{\otimes} \ar[drr]_{\sim}
&{ }\ar@{=>}[ur]_{\lambda}\ar@{=>}[dr]^{\rho}&
 1 \coprod 1 \ar[u]_{e \coprod id}  \ar[d]^{id \coprod e}
\\
&& 1 \coprod 0 }$$ and 
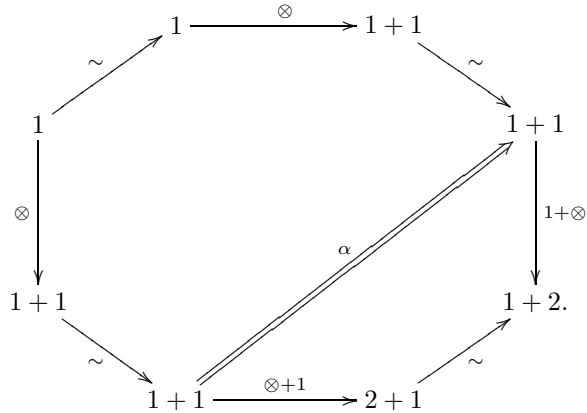
$$\xymatrix{
     &1\ar[rr]^{\otimes} & & 1+1 \ar[dr]^{\sim}
\\
1\ar[ur]^{\sim}\ar[dd]_{\otimes}&  & & &1+1 \ar[dd]^{1+\otimes}
\\
\\
1+1\ar[dr]_{\sim} &&&& 1+2.
\\
& 1+1 \ar[rr]^{\otimes + 1}\ar@{=>}[uuurrr]^{\alpha} && 2+1
\ar[ur]_{\sim} }$$ where the corner isomorphisms $n+m \lra m+n$ is
in $\finbar$. These 2-cells are subject to a unital equation
(left for the reader) and the now-famous pentagon condition:
$$
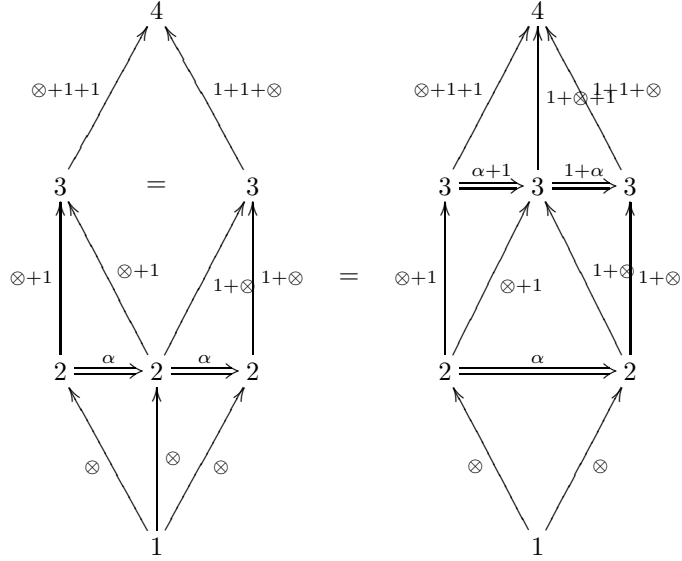
\xymatrix{
      &   4   &      &        &      &   4
\\
\\
3\ar[uur]^{\otimes +1+1}&=       &3\ar[uul]_{1+1+\otimes}  &        &3\ar@{=>}[r]^{\alpha+1}
\ar[ruu]^{\otimes+1+1}  &
3\ar@{=>}[r]^{1+\alpha}\ar[uu]_{1+\otimes+1} &3\ar[uul]_{1+1+\otimes}
\\    &       &      &   =    &
\\
2\ar[uu]^{\otimes+1}\ar@{=>}[r]^{\alpha}&2\ar@{=>}[r]^{\alpha}\ar[uur]_{1+\otimes}\ar[uul]_{\otimes +1}
 &2\ar[uu]_{1+\otimes}
  &&2\ar@{=>}[rr]^{\alpha}\ar[uu]^{\otimes+1}\ar[uur]_{\otimes+1}  &       & 2\ar[uul]_{1+\otimes}
\ar[uu]_{1+\otimes}
\\
\\
      &   1\ar[luu]^{\otimes} \ar[uu]_{\otimes} \ar[ruu]_{\otimes}   &      &
   &      &   1 \ar[luu]^{\otimes}
 \ar[ruu]_{\otimes}
}$$
(We leave out the corner isomorphisms in order to make the
diagram easier to read. However they are important and must be
placed in the definition).
\end{examp}

\begin{examp}
The theory of braided tensor categories $\Tbraid$ and balanced
tensor categories,$\Tbal$,  are easily described in a similar
manner \cite{Joy&Str}.
\end{examp}

\begin{examp}
Associative categories \cite{Paper1} which are monoidal categories
in which the pentagon coherence does not necessarily hold are
described by $\Tassoc$. Similarly, commutative categories \cite{Paper2}
which are braided tensor categories that do not necessarily satisfy the
hexagon coherence condition are described by $\Tcomm$.
\end{examp}

\begin{examp}
Whenever we have a theory with strict associativity, we denote it with
a small ``s'' followed by the usual name e.g. $\Tsmon$, $\Tsbraid$,
$\Tsbal$ etc.
\end{examp}

\begin{defi}
A {\bf 2-theory-morphism} from $\TT 1 $ to $\TT 2 $ is a 2-functor
$G:\TT 1 \lra \TT 2 $ such that $G \circ G_{\TT 1} = G_{\TT 2}$. A
{\bf 2-theory-natural transformation} is a natural transformation
$\gamma :G_1 \Longrightarrow G_2$ between two 2-theory-morphisms.
A {\bf 2-theory-modification} is a modification between two
2-theory-natural transformations.
\end{defi}

We shall denote the 3-category of 2-theories, 2-theory-morphisms,
2-theory-natural transformations and 2-theory-modifications as
$\TTh$.

Here is a diagram of some of the 2-theories and 2-theory
morphisms that we will work with. $$ 
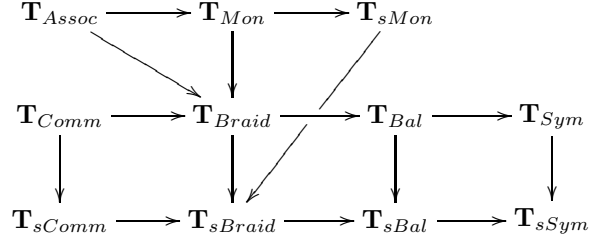
\xymatrix{ \Tassoc \ar[r]
\ar[dr] & \Tmon \ar[r] \ar[d] & \Tsmon \ar[ddl]|\hole
\\
\Tcomm \ar[r] \ar[d]& \Tbraid \ar[r]\ar[d] & \Tbal \ar[r]\ar[d] & \Tsym \ar[d]
\\
\Tscomm \ar[r] & \Tsbraid \ar[r] & \Tsbal \ar[r] & \Tssym
}$$

Many examples of 2-theories  and their morphisms come from one
dimensional theories in the following way. Let $\th $ denote the
usual \cite{Lawvere} 2-category of theories, theory-morphisms and
theory-natural transformations. One can think of $\th $ as a
3-category $\Th $ with only trivial 3-cells. Analogous to the
relationship between sets and topological spaces, we have the
following adjunctions: \be \xymatrix{\label{XY:adj} \Th
\ar@/_1pc/@<-3ex>[rrrr]_c \ar@<2ex>[rrrr]^d_{\bot} &&&&\TTh .
\ar@/_1pc/@<-3ex>[llll]_{\pi_0}^{\bot} \ar@<2ex>[llll]_U^{\bot}
}\ee

$c(T)$ is the 2-theory with the same 1-cells as $T$ and a unique
2-cell between nontrivial 1-cells. $d(T)$ has the same 1-cells as
$T$ and only trivial 2-cells. $U( \T)$ forgets the 2-cells of
$\T$. $\pi_0( \T)$ is a quotient theory of $\T$ where two 1-cells
are set equal if there is a 2-cell between them. These functors
extend in an obvious way to 3-functors. By adjunction we mean a
strict 3-adjunction; that is the universal property is satisfied
by a strict 2-category isomorphism. For example the following
2-categories
 are isomorphic
$$Hom_{\Th}(T, U(\T)) \cong Hom_{\TTh}(d(T), \T)$$

\begin{examp}
$\finbar = d(\fin )$, that is, the theory of categories is
the discrete theory of sets.
\end{examp}

\begin{examp}
$\Tbin = d(T_{Magmas})$.
\end{examp}

\begin{examp}
$d(T_{Monoids})$ is the theory of strict monoidal
categories, $\Tsmon$.
\end{examp}

\begin{examp}
Let $T_{Magmas \bullet}$ be the theory of pointed magmas i.e. the
theory of magmas with a distinguished element. $c(T_{Magmas
\bullet})$ is the 2-theory of symmetric (monoidal ) tensor
categories. Warning: not all operations are supposed to be
isomorphic to one another. In particular, the projections
(inclusions) live in $\finbar$ and are not isomorphic.
\end{examp}

\begin{examp}Let $\Tbraid $ denote the 2-theory of braided tensor categories.
$\pi_0(\Tbraid )$ is
the theory of commutative monoids.
\end{examp}

The units and counits of these adjunctions are of interest.
$\varepsilon : \pi_0 d T \longrightarrow T$,
$\mu : T \longrightarrow U d T$ and
 $ \varepsilon : U c T \longrightarrow T$ are all identity theory-morphisms.
More importantly, $\mu : \T \longrightarrow d \pi_0 \T$ is the
2-theory-morphism corresponding to ``strictification''. Every
2-cell becomes the identity. ``Strictification'' is often used in
coherence theory. Similarly,
 $\mu : \T \longrightarrow c U \T$
might be called ``coherification'' where a 2-theory is forced to
be coherent. $\varepsilon: d U \T \longrightarrow \T$ is the
injection of the 1-theory into the 2-theory.

\begin{defi}
Given a 2-theory $\T$ and a 2-category $\C $ with a product
structure, an ${\bf algebra}$ of $\T$ in $\C $ is a product
preserving 2-functor $F:\T^{op} \longrightarrow \C. $
\end{defi}

\begin{defi}
A {\bf quasi-natural transformation} \cite{Bunge1, Gray}
 $\sigma$ from an algebra $F$ to an to an algebra $F'$ is
\begin{itemize}
\item A family of 1-cells in $\C $, $\sigma_n:F(n) \longrightarrow
F'(n)$ indexed by 0-cells of $\T $. This family must preserve
products i.e. $\sigma_n = (\sigma_1)^n:F(1)^n \longrightarrow
F'(1)^n$.
\item A family of 2-cells in $\C $, $\sigma_f$, indexed by 1-cells
$f:m\longrightarrow n$ of $\T$. $\sigma_f$ makes the following
diagram commute.
\end{itemize}
$$\xymatrix{ & F(1)^n \ar[r]^{\sigma^n}& F'(1)^n\ar[dr]^{\sim}
\\
F(n)\ar[ur]^{\sim} \ar[dd]_{Ff}& & & F'(n)\ar[dd]_{F'f}&
\\
&&& { }
\\
F(m) \ar[dr]_{\sim} &&& F'(m) &
\\
& F(1)^m \ar[r]^{\sigma^m}\ar@{=>}[rruu]^{\sigma_f}  & F'(1)^m
\ar[ur]_{\sim} && }$$

These morphisms must satisfy the following conditions:
\begin{enumerate}
\item If $f$ is in the image of $G_\T:\finbar \longrightarrow \T$,
then $\sigma_f = id$. That is, the quasi-commutative diagram must
commute strictly. This condition includes $\sigma_{id_n} =
id_{\sigma_n}$.
\item $\sigma$ preserves the (co)product structure: $\sigma_{f + f'} =
\sigma_f \times \sigma_{f'}$.  See \cite{Paper3} for an exact
diagram.
\item $\sigma_{g\circ f} = \sigma_f \circ_v \sigma_g $where $\circ_v=\circ_1$ is the
vertical composition of 2-cells.
\item $\sigma$ behaves well with respect to 2-cells of $\T$. See \cite{Paper3} for an exact diagram.
\end{enumerate}
\end{defi}

We shall call $\sigma$ an {\bf iso-}quasi-natural transformation
if $\sigma_f$ is an iso-2-cell for all $f \in_1 \T$.

\begin{defi}
Given two quasi-natural transformations $\sigma , \sigma' : F
\longrightarrow F'$, a ${\bf modification}$ $\Sigma:\sigma
\leadsto \sigma'$ from $\sigma$ to $\sigma'$ is a family of
2-cells $\Sigma_n:\sigma_n \Longrightarrow \sigma'_n$ indexed by
the 0-cells of $\T$. These 2-cells must satisfy the following
conditions:
\begin{enumerate}
\item $\Sigma$ preserves products i.e. $\Sigma_n = (\Sigma_1)^n:(\sigma_1)^n
\Longrightarrow (\sigma_1')^n$.
\item $\Sigma$ behaves well with respect to the 2-cells of $\T $. That is,
if we have
$$\xymatrix{
m\rrtwocell^f_{f'}{\alpha}&&n
} $$
then we have the following ``cube relation'':
\end{enumerate}
\end{defi}
$$ \xymatrix{ F(n) \ar[rr]^{\sigma'_n} \ar[dd]_{id} && F'(n)
\ar[dr]^{F'(f')} \ar[dd]^{id}
 &&
F(n) \ar[rr]^{\sigma'_n} \ar[dd]_{id} \ar[dr]^{F(f')} && F'(n) \ar[dr]^{F'(f')}
\\
 & &{ } & F'(m) \ar[dd]^{id}
& &        F(m) \ar[rr]_{\sigma'_m} \ar[dd]^{id} &{ }\ar@{=>}[u]_{\sigma'_{f'}} & F'(m) \ar[dd]^{id}
\\
F(n) \ar[rr]_{\sigma_n} \ar[rd]_{F(f)} &{ }\ar@{=>}[ur]^{\Sigma_n}& F'(n) \ar[rd]^{F'(f)}\ar@{=>}[r]^{F'(\alpha)}
&{ }\ar@{}[r]^= & F(n) \ar[dr]_{F(f)}\ar@{=>}[r]^{F(\alpha)}&&&{ }
\\
& F(m) \ar[rr]_{\sigma_m} &{ }\ar@{=>}[u]^{\sigma_f} &F'(m) &&
F(m) \ar[rr]_{\sigma_m} &{ }\ar@{=>}[ur]^{\Sigma_m}&F'(m) }$$

For a given 2-theory $\T$ and a 2-category $\C$ with a product
structure, we denote the 2-category of algebras, quasi-natural
transformations and modifications as ${\bf 2Alg(T, C)}$. We shall
denote the locally full sub-2-category of algebras, {\bf
iso-}quasi-natural transformations and modifications as ${\bf
2Alg^i(T, C)}$.

(A quasi-natural transformation is a way of having an operations
preserved up to a 2-cell. Many times in coherence theory, one
wants some operations to be preserved up to a 2-cell and some
operations to be preserved strictly. This is done in
\cite{Paper3} with the notion of a relative quasi-natural
transformation. We demand two 2-theories $\TT 1, \TT 2$ and a
2-theory-morphism between them $G:\TT 1 \lra \TT 2$. $\TT 1$
controls which operations in $\TT 2$ should be preserved strictly.
A relative quasi natural transformation between two $\TT 2$
algebras is a quasi-natural transformation $\sigma : F \Lra F'$
such that $\sigma \circ G^{op}$ is a natural transformation (not
quasi) from $F \circ G^{op}$ to  $F' \circ G^{op}$. $$ \xymatrix{
\T_1^{op}\ar[rr]^{G^{op}} && \T_2^{op}\rrtwocell^F_{F'}{\sigma} &&
\C }.$$ ${\bf 2Alg_G(\TT 2, \C)}$ has the same algebras as ${\bf
2Alg(T_2, C)}$ but with only relative quasi-natural
transformations between them. The present paper will not use this
notion. However, we will discuss the relationship between relative
quasi-natural transformations and relative homotopy theory in
Section 6.)

Fixing $\C = \Cat$ we can extend ${\bf 2Alg( - , \Cat)}$ to be a
3-functor from $\TTh^{op}$ to $\TCatcat$ in the obvious way.

It is common to look at the algebras of one theory in the
category of algebras of another theory. The theory of such
algebras is given as the Kronecker product of the two theories.

The Kronecker product of (1-)theories is a well understood
coherent symmetric monoidal 2-bifunctor $\otimes_K:\th \times \th
\longrightarrow \th$. Let $T_1$ and $T_2$ be two theories. $T_1
\otimes_K T_2$ is a theory that satisfies the universal property
$$ Alg(T_1 \otimes_K T_2, C) \cong  Alg(T_1, Alg(T_2,C)).$$ $T_1
\otimes_K T_2$ is constructed as follows. Construct the coproduct
$T_1 \coprod T_2$ in the category of theories (pushout in
$\Cat$.) Place a congruence on  $ T_1 \coprod T_2$ such that for
all
 $f:m\longrightarrow m'$ in $T_1$ and $g:n \longrightarrow n'$  in
$T_2$ the diagram $$ \xymatrix{
        &n^m \ar[rr]^{g^m} &&n'^m \ar[dr]^{\sim}
\\
m^n \ar[dd]_{f^n}\ar[ur]^{\sim} &&&&m^{n'} \ar[dd]^{f^{n'}}
\\
\\
{m'}^n \ar[dr]^{\sim}&&&&{m'}^{n'}
\\
        &n^{m'} \ar[rr]^{g^{m'}} &&{n'}^{m'} \ar[ur]^{\sim}
}$$ commutes. We have a  full theory-morphism $ T_1 \coprod T_2
\longrightarrow T_1 \otimes_K T_2.$

We will work with a two-dimensional analogue to the Kronecker
product. (See \cite{Paper3} for more details.)

\begin{defi}
A {\bf (2-)Kronecker product} of 2-theories is a 3-bifunctor $$\KP
: \TTh \times \TTh \longrightarrow \TTh$$ that satisfies the
following universal property: for all $$\xymatrix{\TT 1 & &
\finbar \ar[ll]_{G_{\TT 1}} \ar[rr]^{G_{\TT 2}} && \TT 2}$$ there
is an induced $$\xymatrix{ \finbar \ar[rr]^{G_{\T_1}}
\ar[dd]_{G_{\T_2}} \ar[ddrr]^{G_{\T_1} \KP G_{\T_2}} &&\T_1
\ar[dd]
\\
\\
\T_2 \ar[rr]&& \T_1 \KP \T_2  }$$ and for all 2-categories with
finite products $\C$, an isomorphism  of 2-categories $$ {\bf
2Alg}(\TT 1 \KP \TT 2, \C) \cong {\bf 2Alg}(\TT 1, {\bf 2Alg}(\TT
2, \C))$$ which is natural for all cells in $ \TTh$ and $\C$.
\end{defi}

It will be helpful to examine the naturality conditions in terms
of 1-cells of $\TTh$. Let $F_1:\TT 1 \lra \TT 1'$ and $F_2: \TT 2
\lra \TT 2 '$ be 2-theory morphisms. By the functoriality of
$\KP$ there is an induced 2-theory-morphism $F_1 \KP F_2: \TT 1
\KP \TT 2 \lra \TT 1' \KP \TT 2'$ such that $${\bf  \xymatrix{
2Alg(\TT 1' \KP \TT 2' ,\C) \ar[rr]^\cong \ar[dd]_{2Alg(F_1 \KP
F_2, \C)}&& 2Alg(\TT 1',
2Alg(\TT 2', \C)) \ar[dd]^{2Alg(F_1,2Alg(F_2, \C))} \\
\\
2Alg(\TT 1 \KP \TT 2 ,\C) \ar[rr]^\cong && 2Alg(\TT 1, 2Alg(\TT
2, \C))
      } }$$

In order to construct $\TT 1 \KP \TT 2$, we take the coproduct
 $\TT 1 \coprod_{\finbar} \TT 2$ in $\TTh$ and
we freely add in the following 2-cells: For every $f:m'
\longrightarrow m$ in $\TT 1$ and $g:n' \longrightarrow n$ in
$\TT 2$  we add the 2-cell $\delta_{\T_1, \T_2} (f,g)$ that makes
the following diagram commute: \be  \xymatrix{
        &n^m \ar[rr]^{g^m} &&n'^m \ar[dr]^{\sim}
\\
m^n \ar[dd]_{f^n}\ar[ur]^{\sim} &&&&m^{n'} \ar[dd]^{f^{n'}}
\\
&&&&{ }
\\
{m'}^n \ar[dr]^{\sim}&&&&{m'}^{n'}
\\
        &n^{m'} \ar[rr]^{g^{m'}} \ar@{=>}[rrruuu]^{\delta_{\T_1, \T_2}
(f,g)}&{ }&{n'}^{m'} \ar[ur]^{\sim} }\ee

The $\delta$`s must satisfy the following coherence conditions
that are compatible to the four coherence conditions in the
definition of a quasi-natural transformation.
\begin{enumerate}
\item If $f$ is in the image of $G_1$, then $\delta (f,g)$ must be set to the identity.
\item $\delta$ must preserve products in $f$.
\item  $\delta (f\circ f',g) = \delta(f,g) \circ_v \delta(f',g) $
\item $\delta$ must preserve 2-cells.
\end{enumerate}

The fact that there is choice in the construction of $\TT 1 \KP
\TT 2$, should not disturb the reader since we never claimed that
$\TT 1 \KP \TT 2$ should be unique. Rather, it should be unique up
to a (2-)isomorphism. In order to see that our construction of
$\TT 1 \KP \TT 2$ satisfies the universal properties demanded of
it, we must realize that our construction was made to mimic the
definition of a quasi-natural transformation in our 2-categories
of algebras.

Given $F_1:\TT 1 \lra \TT 1'$ and $F_2: \TT 2 \lra \TT 2 '$, we
construct $F_1 \KP F_2$ as follows: construct $F_1 \coprod F_2:
\TT 1 \coprod \TT 2 \lra \TT 1' \coprod \TT 2'$ and define $F_1
\KP F_2$ on $\delta_{\TT 1, \TT 2}(f,g)$ as

\be (F_1 \KP F_2)(\delta_{\T_1, \T_2}(f,g)) \quad = \quad
\delta_{\T_1', \T_2'}(F_1(f), F_2(g)).\ee

\section{Biequivalences}
The notion of a biequivalence is central to this paper. A
biequivalence is a 2-categorical generalization of an
equivalence. We know of no original sources for the idea or the
name. The properties of  biequivalences seem to be well known
folklore that remains unwritten. We shall be explicit with some
of these properties.

\begin{defi}
Within a 2-category $\A$, a 1-cell $f: a \lra a'$ is an {\bf
equivalence} if there exists a 1-cell $g:a' \lra a$ and
isomorphic 2-cells $\eta:id_a \lrat g\circ f$ and $\varepsilon: f
\circ g \lrat id_{a'}$. We denote such an equivalence as
$(f,g,\eta,\varepsilon):a \lra a'$.
\end{defi}

\begin{defi}
Within a 3-category $\tilde{\C}$, a 1-cell $F:\A
\lra \B$ is a {\bf biequivalence} if there is a 1-cell $G:\B \lra
\A$ with an equivalence $(\eta, \delta, \Phi, \Xi):id_\A \lra
G\circ F$ in $\tilde{\C}(\A, \A)$ and an equivalence
$(\varepsilon, \zeta, \Psi, \Omega):F \circ G \lra id_\B$ in
$\tilde{\C}(\B, \B)$.
\end{defi}

Let us look at this definition in more detail. $F:\A \lra \B$ is
a biequivalence if $F$ is part of a 10-tuple $$(F,G, \eta, \delta,
\varepsilon, \zeta, \Phi, \Xi, \Psi, \Omega):\A \lra \B$$ where

$$F:\A \lra \B \qquad \qquad \qquad \qquad  G:\B\lra \A$$

$$\eta:id_\A \lra G \circ F \qquad \delta: G\circ F \lra id_\A
\qquad \varepsilon:F\circ G \lra id_\B \qquad \zeta:id_\B \lra
F\circ G$$

$$\Phi:id_\A \lrat \delta \circ \eta  \quad \Xi: \eta\circ \delta
\lrat id_{G \circ F} \qquad \qquad  \Psi:id_{F\circ G} \lrat \zeta
\circ \varepsilon \quad \Omega:\varepsilon \circ \zeta \lrat
id_\B$$

(Ross Street \cite{Street} has written a general definition of a
$k$-equivalence between two $n$-categories where $k \le n+1$.
Using that language, a biequivalence is a 3-equivalence of two
2-categories. See \cite{Paper4} for general properties of a
$k$-equivalence.)

\begin{prop}
Let $(F,G, \eta, \delta, \varepsilon, \zeta, \Phi, \Xi, \Psi,
\Omega):\A \lra \B$ be a biequivalence. The following are also
biequivalences.

1) $(F,G, \delta,\eta,  \varepsilon, \zeta, \Xi^{-1}, \Phi^{-1},
\Psi, \Omega):\A \lra \B$

2) $(F,G, \eta, \delta, \zeta, \varepsilon,  \Phi, \Xi,
\Omega^{-1},\Psi^{-1}):\A \lra \B$

3) $(G,F, \zeta,\varepsilon,\delta, \eta, \Psi,\Omega, \Phi, \Xi,
):\B \lra \A$
\end{prop}

\begin{prop}
Let $(F,G, \eta, \delta, \varepsilon, \zeta, \Phi, \Xi, \Psi,
\Omega):\A \lra \B$ and $(F',G', \eta', \delta', \varepsilon',
\zeta', \Phi', \Xi', \Psi', \Omega'):\B \lra \C$ be
biequivalences. Then $(F' \circ F,\quad G\circ G',\quad (G\eta'
F)\circ \eta,\quad \delta\circ (G\delta' F),\quad \varepsilon'
\circ (F'\varepsilon G'),\quad (F' \zeta G')\circ \zeta',\quad
\Phi'',\quad \Xi'',\quad \Psi'',\quad \Omega''):\A \lra \C$ is
also a biequivalence where

$\Phi''=[(\delta  G \delta')\zeta(\eta' F \eta)] \circ[(\delta G
)\Phi'(F \eta)] \circ[(\delta)\eta(\eta)] \circ \Phi$

$\Xi''=\delta \circ[(G)\Xi'(F)]\circ[(G \eta' )\varepsilon (
\delta ' F)]\circ[(G \eta' F)\Xi(G \delta ' F)]$

$\Psi''=[(F'\zeta G') \Psi'(F'\varepsilon G')]\circ[(F'\zeta )
\eta(\varepsilon G')]\circ[(F')\Psi( G')]\circ\varepsilon'$

$\Omega''=\Omega' \circ
[(\varepsilon')\varepsilon(\zeta')]\circ[(\varepsilon' F)\Omega(
G' \zeta')]\circ[(\varepsilon' F \varepsilon)\delta' (\zeta G'
\zeta' )].$
\end{prop}

\begin{prop}
Let $(F_1,G_1, \eta_1, \delta_1, \varepsilon_1, \zeta_1, \Phi_1,
\Xi_1, \Psi_1, \Omega_1):\A_1 \lra \B_1$ and $$(F_2,G_2, \eta_2,
\delta_2, \varepsilon_2, \zeta_2, \Phi_2, \Xi_2, \Psi_2,
\Omega_2):\A_2 \lra \B_2$$ be biequivalences. Then $$((F_1 \coprod
F_2) ,(G_1 \coprod G_2) , (\eta_1 \coprod \eta_2) , (\delta_1
\coprod \delta_2) , (\varepsilon_1 \coprod \varepsilon_2),
(\zeta_1 \coprod \zeta_2) ,$$ $$ (\Phi_1 \coprod \Phi_2) , (\Xi_1
\coprod \Xi_2) , (\Psi_1 \coprod \Psi_2) , (\Omega_1 \coprod
\Omega_2)):(\A_1 \coprod \A_2) \lra (\B_1 \coprod \B_2)$$ is also
a biequivalence. (Similarly for products.)
\end{prop}

\begin{prop}
Let $F:\A \lra \B$ be a biequivalence and $G:\B \lra \A$ be a
2-functor such that

1)$G \circ F = id_\A$ (i.e. $F$ is the inclusion of a full
sub-2-category) then there is a biequivalence $(F,G, id, id,
\varepsilon, \zeta, id,id, \Psi, \Omega):\A \lra \B$

2)$F \circ G = id_\B$ (i.e. $F$ is a surjection of 2-categories)
then there is a biequivalence $(F,G, \eta, \delta, id, id, \zeta,
\Phi, \Xi, id, id):\A \lra \B$

3)$G \circ F \cong id_\A$ (i.e. $F$ is ``almost'' the inclusion
of a full sub-2-category) then there is a biequivalence $(F,G,
\eta, \eta^{-1}, \varepsilon, \zeta, id, id, \Psi, \Omega):\A \lra
\B$

4)$G \circ F \cong id_\B$ (i.e. $F$ is ``almost'' a surjection of
2-categories) then there is a biequivalence $(F,G, \eta, \delta,
\varepsilon, \varepsilon^{-1}, \Phi, \Xi, id, id):\A \lra \B$
\end{prop}

\begin{defi}
A biequivalence is called an {\bf adjoint biequivalence} if $$ F
\eta \cong \zeta F \qquad F \delta \cong \varepsilon F \qquad G
\zeta \cong \eta G \qquad G \varepsilon \cong \delta G$$
\end{defi}

\begin{prop}
If $(F,G, \eta, \delta, \varepsilon, \zeta, \Phi, \Xi, \Psi,
\Omega):\A \lra \B$ is a biequivalence then the following are
adjoint biequivalences

1) $(F,G, \eta', \delta', \varepsilon, \zeta, \Phi', \Xi', \Psi,
\Omega):\A \lra \B$ where

$\qquad \eta'_a = \delta_{GFa} \circ G \zeta_{Fa} \circ \eta_a$

$\qquad \delta'_a = \delta_{a} \circ G \varepsilon_{Fa} \circ
\eta_{GFa}$

$\qquad \Phi'_a = [\delta_a \circ G \varepsilon_{Fa} \circ
\Xi^{-1}_{GFa} \circ G \zeta_{Fa} \eta_a] \circ [\delta_a \circ G
\Omega^{-1}_{Fa} \eta_a]\circ \Phi_a$

$\qquad \Xi'_a = \Phi^{-1}_{GFa} \circ [\delta_{GFa} \circ G
\Psi^{-1}_{Fa} \circ \eta_{GFa}] \circ [\delta_{GFa} \circ G
\zeta_{Fa} \circ \Xi_a \circ G\varepsilon_{Fa} \circ \eta_{GFa}]$

2) $(F,G, \eta, \delta, \varepsilon', \zeta', \Phi, \Xi, \Psi',
\Omega'):\A \lra \B$ where

$\qquad \varepsilon'_b = \varepsilon_b \circ F \delta_{Gb} \circ
\zeta_{FGb}$

$\qquad \zeta'_b = \varepsilon_{FGb} \circ F \eta_{Gb} \circ
\zeta_{b}$

$\qquad \Phi'_b = [\varepsilon_{FGb} \circ F \eta_{Gb} \circ
\Psi^{-1}_{b} \circ F \delta_{Gb} \zeta_{FGb}] \circ
[\varepsilon_{FGb} \circ F \Xi^{-1}_{Gb} \zeta_{FGb}]\circ
\Omega^{-1}_{FGb}$

$\qquad \Omega'_b = \Omega_{b} \circ [\varepsilon_{b} \circ F
\Phi^{-1}_{Gb} \circ \zeta_{b}] \circ [\varepsilon_{b} \circ F
\delta_{Gb} \circ \Psi^{-1}_{FGb} \circ F\eta_{Gb} \circ
\zeta_{b}]$

In other words, Any biequivalence can be made into an adjoint
biequivalence.
\end{prop}

{\bf Proof.} We shall focus on part 1). Part 2) is
completely symmetric. The following diagram will help the reader
see $\Phi'$ more clearly.

$$
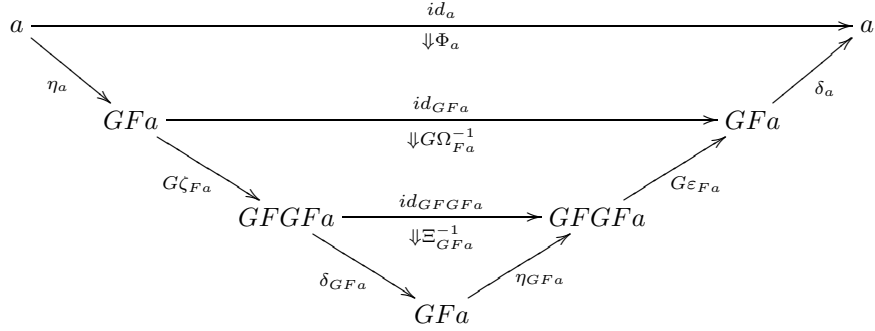
\xymatrix{a \ar[rrrrrr]^{id_a}_{\Downarrow \Phi_a}
 \ar[rd]_{\eta_a} & & & & & & a
\\
& GFa \ar[rrrr]^{id_{GFa}}_{\Downarrow G\Omega^{-1}_{Fa}}
\ar[rd]_{G\zeta_{Fa}} &&&& GFa \ar[ru]_{\delta_a}
\\
&&GFGFa \ar[rr]^{id_{GFGFa}}_{\Downarrow \Xi^{-1}_{GFa}}
\ar[rd]_{\delta_{GFa}} && GFGFa \ar[ru]_{G\varepsilon_{Fa}}
\\
&&&GFa \ar[ru]_{\eta_{GFa}}
 }$$

Similar diagrams can be drawn for $\Xi', \Psi'$ and $\Omega'$. In
order to see that this new biequivalence satisfies the first
adjointness condition ($F\eta' \cong \zeta_F$), consider the
following diagram $$
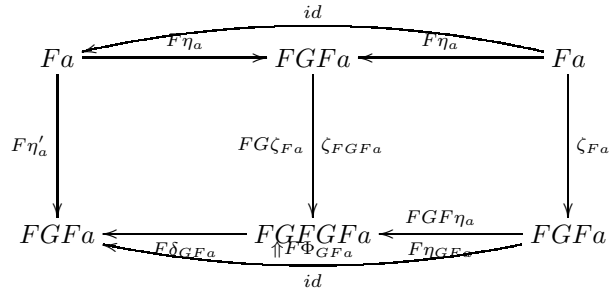
\xymatrix{Fa \ar[dd]_{F\eta'_a}
\ar[rr]^{F\eta_{a}}&& FGFa\ar[dd]_{FG\zeta_{Fa}}^{\zeta_{FGFa}} &&
Fa \ar@/_1pc/[llll]_{id} \ar[dd]^{\zeta_{Fa}}\ar[ll]_{F\eta_{a}}
\\
\\
FGFa && FGFGFa \ar[ll]^{ F\delta_{GFa}}&& FGFa
\ar@/^1pc/[llll]^{id}_{\Uparrow F\Phi_{GFa}}
\ar[ll]^{F\eta_{GFa}}_{FGF\eta_{a}}}$$

Where the left square commutes by the definition of $F\eta'_a$
and the right square commutes by the naturality of $\zeta$. In
detail:
\begin{eqnarray*}
F\eta'_a & = & F\delta_{GFa} \circ FG \zeta_{Fa} \circ F \eta_a
\quad \mbox{  by definition of $\eta'$} \\
& = & F\delta_{GFa} \circ F \eta_{GFa} \circ \eta_{Fa}
\quad \mbox{  by naturality of $\zeta$} \\
&\cong& id \circ \zeta_{Fa} \qquad \mbox{       by $F
\Phi^{-1}_{GFa}$}\\
& = & \zeta_{Fa}.\end{eqnarray*}

The following three diagrams show the other three adjointness
conditions:
 $$\xymatrix{FGFa
\ar[dd]_{F\delta'_a}
\ar@/^1pc/[rrrr]^{id}_{\Downarrow F \Phi_{GFa}} 
\ar[rr]^{F\eta_{GFa}}&&
FGFGFa\ar[dd]_{FG\varepsilon_{Fa}}^{\varepsilon_{FGFa}}
\ar[rr]^{F\delta_{GFa}}_{FGF\delta_{a}}&& FGFa
\ar[dd]^{\varepsilon_{Fa}}
\\
\\
Fa 
&& FGFa \ar[ll]^{F\delta_{a}} \ar[rr]_{F\delta_{a}}&& Fa }$$

 $$\xymatrix{Gb\ar[dd]_{\eta'_{Gb}}
\ar@/^1pc/[rrrr]^{id}_{\Downarrow \Phi_{Gb}}
\ar[rr]^{\eta_{Gb}}&&
 GFGb\ar[dd]_{G\zeta_{FGb}}^{GFG\zeta_b}\ar[rr]^{\delta_{Gb}}
&& Gb \ar[dd]^{G \zeta_b}
\\
\\
GFGb 
&& GFGFGb \ar[ll]^{
GF\delta_{Gb}}\ar[rr]^{\delta_{GFGb}}_{GF\delta_{Gb}}&& GFGb }$$

$$\xymatrix{GFGb \ar[dd]_{\delta'_{Gb}}\ar[rr]^{GF\eta_{Gb}}&&
GFGFGb\ar[dd]_{G\varepsilon_{FGb}}^{GFG\varepsilon_b} && GFGb
\ar[dd]^{G\varepsilon_b}
\ar[ll]^{\eta_{GFGb}}_{GF\eta_{Gb}}
\\
\\
Gb && GFGb \ar[ll]^{ \delta_{Gb}}&&
Gb
\ar@/^1pc/[llll]^{id}_{\Uparrow \Phi_{Gb}} 
\ar[ll]^{\eta_{Gb}}
}$$ $\Box$
\begin{prop}
Let $\A$ and $\B$ be 2-categories and $F: \A \lra \B$ be a
2-functor. F is a biequivalence iff

1) for all $b \in_0 \B$, there exists an $a \in_0 \A$ such that
$F(a)$ is equivalent to $b$

and

2) for all $a,a' \in_0 \A$, $F:\A(a,a') \lra \B(F(a), F(a')) $ is
an equivalence.
\end{prop}
{\bf Proof.} Assume $F$ is part of a an adjoint biequivalence then
for all $b \in \B$, $Gb \in_0 \A$ and we have the equivalence
$(\varepsilon_b, \zeta_b, \Psi_b, \Omega_b):FG(b) \lra b$. The
equivalence for part $2)$ is given as $$ (F, (\delta_{a'} \circ
(-) \circ \eta_a) \circ G, \Phi_{a'}, \Omega_{Fa'}^{-1}):\A(a,a')
\lra \B(F(a),F(a')).$$

In order to understand the unit of this equivalence, consider the
following diagram: $$\xymatrix{ a \ar[rr]^{\eta_a} \ar[dd]_f
\ar@/^1pc/[rrrr]^{id}_{\Downarrow \Phi_a} && GFa
\ar[rr]^{\delta_a}
\ar[dd]_{GFf} && a \ar[dd]^f \\ \\
a' \ar[rr]_{\eta_{a'}} \ar@/_1pc/[rrrr]_{id}^{\Uparrow \Phi_{a'}}
&& GFa' \ar[rr]_{\delta_{a'}} && a' }$$ Let $f\in \A(a,a')$, then

\begin{eqnarray*}
f & = & id_{a'} \circ f \\ & \cong & \delta_{a'} \circ \eta_{a'}
\circ f \quad \mbox{ by } \Phi_{a'} \\
 & = & \delta_{a'} \circ GFf \circ \eta_a \quad \mbox{ by naturality
of }\eta \\
 & = & [(\delta_{a'} \circ (-) \circ \eta_a) \circ G \circ F] f
 \quad
\mbox{ in } \A(a,a').\end{eqnarray*}

In order to understand the counit of this equivalence, consider
the following diagram: $$\xymatrix{ Fa \ar[rr]^{\zeta_{Fa}}
\ar[dd]_g \ar@/^1pc/[rrrr]^{id_F}_{\Uparrow \Omega_{Fa}} && FGFa
\ar[rr]^{\varepsilon_{Fa}}
\ar[dd]_{FGg} && Fa \ar[dd]^g \\ \\
Fa' \ar[rr]_{\zeta_{Fa'}} \ar@/_1pc/[rrrr]_{id}^{\Downarrow
\Omega_{Fa'}} && FGFa' \ar[rr]_{\varepsilon_{Fa'}} && Fa' }$$ Let
$g\in \B(Fa,Fa')$, then
\begin{eqnarray*}
g & = & id_{Fa'} \circ g \\ & \cong & \varepsilon_{Fa'} \circ
\zeta_{Fa'} \circ g \quad \mbox{ by } \Omega_{Fa'}^{-1} \\ & = &
\varepsilon_{Fa'} \circ FGg \circ \zeta_{Fa} \quad \mbox{ by
naturality of }\zeta \\ & \cong & F \delta_{a'} \circ FGg \circ
F\eta_a \quad \mbox{ by adjoint biequivalence }\\ & = & F\circ
(\delta_{a'} \circ G(-) \circ \eta_a) (g) \quad \mbox{ in }
\B(Fa,Fa').
\end{eqnarray*}

Conversely, assume that for all $b \in_0 \B$ there is an $a_b
\in_0 \A$ such that $(f_b, G_b, \eta_b, \varepsilon_b):F(a_b)
\lra b$ is an equivalence and for all $a,a'\in_0 \A$ there is an
equivalence
$$(f_{a,a'},g_{a,a'},\eta_{a,a'},\varepsilon_{a,a'}):\A(a,a') \lra
\B(Fa,Fa').$$ Then we shall define $G:\B \lra \A$ as follows

$G(b) \quad = \quad a_b$

$G(h:b \lra b') \quad = \quad g_{a_b,a_{b'}}(g_{b'} \circ h \circ
f_b)$

$G(\beta:h \lra h') \quad = \quad g_{a_b,a_{b'}}(g_{b'} \circ
\beta \circ f_b)$

Define $\alpha: id_\A \lra G \circ F$ and $\gamma: F \circ G \lra
id_\B$ as $$\alpha_a = g_{a,GFa}(g_{Fa}) \qquad \gamma_b =
g_{a,GFa}(f_{Fa}).$$
\begin{eqnarray*}
(\gamma \circ \alpha)_a & = & \gamma_a \circ \alpha_a \\ & = &
g_{a,GFa}(f_{Fa}) \circ g_{a,GFa}(g_{Fa})\\ & =&
g_{a,GFa}(f_{Fa}\circ g_{Ga})\\ & \cong & g_{a,GFa}(id_{Fa}) \quad
\mbox{ by } \varepsilon_{Fa}\\ & =& id.
\end{eqnarray*}
Similarly for $\alpha \circ \gamma$. Hence we have the
sought-after biequivalence. $\Box$

Let us restrict our attention from general 2-functors to
2-theory-morphisms. Since 2-theory-morphisms are bijective on
0-cells, by Proposition 6 we have that $F:\T_1 \lra \T_2$ is a
biequivalence iff $F:\T_1(m,n) \lra \T_2(m,n)$ is an equivalence
of categories for all $m,n \in \N$. A nice example of this is
when $\T_1 = \Tmon$, $\T_2 = \Tsmon$, $m = 1$ and $n = 4$. We
then have the following two equivalent categories.

$$
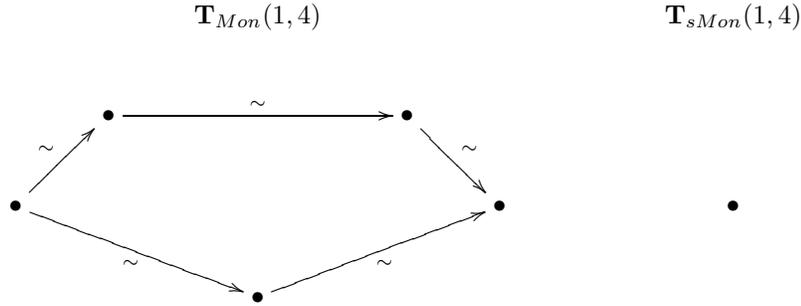
\xymatrix{ && \Tmon(1,4)&&&&\Tsmon(1,4)
\\
& \bullet \ar[rr]^\sim&&\bullet\ar[rd]^\sim \\
\bullet\ar[rrd]_\sim \ar[ru]^\sim &&&& \bullet && \bullet\\
&&\bullet\ar[rru]_\sim}$$

(Notice the importance of the morphisms being isomorphisms. This
tells us that the Mac Lane's coherence theorem would not apply for
categories with a tensor product and a non-isomorphism $A \otimes
(B \otimes C) \lra (A \otimes B) \otimes C$ satisfying the
pentagon condition.)

Similarly, there are obvious biequivalences between the 2-theory
of braided monoidal categories and the 2-theory of strictly
associative braided monoidal categories.

One can go on to prove a more general statement: For any
algebraic 2-theory $\Tx$ with an inclusion of $\Tmon$, we have
the following pushout

$$\xymatrix{ \Tmon \ar@{^{(}->}[rr] \ar[dd] && \Tx \ar[dd]\\
\\
\Tsmon \ar@{^{(}->}[rr]&& \Tsx }$$ where $\Tsx$ is the 2-theory
$\Tx$ with strict associativity. Since the left-hand side of the
pushout is a biequivalence, the right-hand side is also.

The importance of biequivalences for the semantics of coherence
theory is the following proposition.
\begin{prop}
$F:\T_1 \lra \T_2$ is a biequivalence iff $$F^*={\bf 2Alg^i(
}F,{\bf \Cat):2Alg^i(\T_2, \Cat)  \lra 2Alg^i(\T_1, \Cat)}$$ is a
biequivalence.
\end{prop}
{\bf Proof.} If $F$ is a biequivalence, then simply by the
3-functoriality of $F^*={\bf 2Alg^i( - , \Cat)}$ the conclusion
follows.

Conversely, assume $F^*$  is a biequivalence. By Proposition 5,
we may assume that $F^*$ is an adjoint biequivalence. In
\cite{Paper3} Proposition 2, we proved that for any 2-theory $\T$,
$\To(n,-):\T \lra \Cat$ is the free $\T$-algebra on $n$
generators. Since $F^*$ is an adjoint biequivalence, it is not
hard to show that $$F^*(\To_2(n,-)) \cong \To_1(n,-).$$ In order
to show that $F:\T_1 \lra \T_2$ is a biequivalence, it suffices
to show that for all $m,n \in \N$, $F:\T_1(m,n) \lra \T_2(m,n)$
is an equivalence (Proposition 6).
\begin{eqnarray*}
\To_2(m,n) &\simeq& {\bf 2Alg^i(} \T_2,
\Cat)(\To_2(m,-),\To_2(n,-)) \\ & &\qquad \mbox{ by quasi-Yoneda
lemma
\cite{Paper3} Proposition 1 (e)} \\
 & \cong &{\bf 2Alg^i(} \T_1,
\Cat)(F^*(\To_2(m,-)),F^*(\To_2(n,-)))\\& & \qquad \mbox{ since
$F^*$ is
a biequivalence, Proposition 6 } \\
& \cong &{\bf 2Alg^i(} \T_1, \Cat)(\To_1(m,-)),\To_1(n,-))\\& &
\qquad \mbox{ since $F^*$ is
an adjoint biequivalence} \\
&\simeq& \To_1(m,n)\quad \mbox{ by quasi-Yoneda lemma
\cite{Paper3} Proposition 1 (e)}
\end{eqnarray*}
where $\simeq$ means equivalence. $\Box$

Many coherence results simply fall out of Proposition 7. For
example, from the fact that $\Tmon$ is biequivalent to $\Tsmon$
and Propositions 6 and 7 we have that every monoidal category is
tensor equivalent to a strict monoidal category. More exotic
statements can be asserted about higher cells in ${\bf 2Alg(\Tmon,
\Cat)}$.

\section{Quillen Model Category Structure}
In this section, we shall show that the category of $\TTh$ has a
functorial closed Quillen model category (FCQMC) structure. We
must point out that we are (perhaps wrongly) ignoring the higher
categorical structure in the 3-category of $\TTh$. We will only
talk of the (1-)category of 2-theories and 2-theory-morphisms.
There will be more about this omission in Section 6.

A category is given a FCQMC structure by describing three
subclasses of morphisms in the category (weak equivalences,
fibrations and cofibrations) and showing that they satisfy
certain axioms.

\begin{defi} {\bf Weak equivalence} are 2-theory-morphisms that
are biequivalences. {\bf Fibrations} are 2-theory-morphisms
$F:\T_1 \lra \T_2$ satisfying the iso-2-cell lifting property.
That is, for all $f \in_1 \T_1$ and for all iso-2-cells $\beta:Ff
\Lrat g$ in $\T_2$ there is an iso-2-cell $\alpha : f \Lrat f'$
in $\T_1$ such that $F(\alpha)=\beta$. {\bf Cofibrations} are
2-theory-morphisms that are injective on 1-cells. $F:\T_1 \lra
\T_2$ are trivial fibrations (resp. trivial cofibrations) if $F$
is both a fibration (resp. cofibration) and a weak equivalence.
\end{defi}

Since any 2-theory-morphism $F:\T_1 \lra \T_2$ is bijective on
0-cells, one can describe these classes of maps by looking at the
1-functors $\dot{F}:\T_1(m,n) \lra \T_2(m,n)$ for all $m,n \in
\N$. $F$ is a weak equivalence iff the $\dot{F}$'s are
equivalences of categories. $F$ is a fibration iff the
$\dot{F}$'s have the isomorphisms lifting property. $F$ is a
cofibration iff the $\dot{F}$'s are injective on 0-cells.
Furthermore, since a product of 2-theory-morphisms is a weak
equivalence (resp. fibration, cofibration) iff each of its terms
is a biequivalence (resp. fibration, cofibration) and since
$\T(m,n) = \T(1,n)^m$ we need only look at $\dot{F}:\T_1(1,n)
\lra \T_2(1,n)$ for all $n \in \N$ in order to classify $F$. (For
the reader who likes the language of operads, as opposed to
2-theories, we have just reduced the problem of classifying
2-theory-morphisms to be a problem of classifying
2-operad-morphisms. 2-operads are operads in $\Cat$. With an
understanding of this paragraph, one can simply rewrite the
entire paper in the language of operads). This FCQMC structure is
closely associated to Rezk's \cite{Rezk} FCQMC structure on
$\Cat$.

 In order to have a FCQMC structure on the
category $\TTh$, these three subclasses of morphisms must satisfy
the following five axiom schemes.

{\bf Limits and Colimits}. $\TTh$ must have all finite limits and
colimits. These (co)limits are constructed like (co)limits in
$\Cat$ and $\Th$. It is worth pointing out that the initial
2-theory is $\finbar$. The terminal 2-theory, $\TT t$, has $\TT
t(m,n)= \{* \}$ for all $m,n \in \N$. $\TA(\TT t, \Cat)$ has only
one object. The algebras of $\T_1 \coprod \T_2$ are categories
with both $\T_1$ and $\T_2$ structures. The algebras of $\T_1
\times \T_2$ are categories with either a $\T_1$ or a $\T_2$
structure. $\Box$

{\bf Two out of Three}. If $F,G$ or $G \circ F$ are 2-theory morphisms
 and any two of them are biequivalences, then so is the third. If $F$ and $G$ are
 biequivalences then so is $G \circ F$ by Proposition 2. If
$(F,F', \eta_1, \delta_1, \varepsilon_1, \zeta_1, \Phi_1, \Xi_1, \Psi_1,
\Omega_1):\T_1 \lra \T_2$
and
$(G\circ F,(G\circ F)', \eta_2, \delta_2, \varepsilon_2, \zeta_2, \Phi_2, \Xi_2, \Psi_2,
\Omega_2):\T_1 \lra \T_3$ are biequivalences as in

$$\xymatrix{\T_1 \ar[rr]^F \ar@/^1pc/[rrrr]^{G \circ F}&& \T_2
\ar[rr]^G \ar@/^1pc/[ll]_{F'}&& \T_3 \ar@/^1pc/[llll]^{(G \circ
F)'}. }$$ Then we have the following equivalence $$ (\alpha ,
\beta ,\Gamma, \Theta):F' \lra (G \circ F)'G$$ where

$\alpha \quad=\quad [(g \circ F)' G]\varepsilon_2 \quad \circ_1
\quad \eta_1 F'$

$\beta \quad=\quad [(g \circ F)' G]\zeta_2 \quad \circ_1 \quad
\delta_1 F'$

$\Gamma \quad=\quad [(g \circ F)' G]\Psi_2 \quad \circ_1 \quad
\Phi_1 F'$

$\Theta \quad=\quad [(g \circ F)' G]\Omega_2 \quad \circ_1 \quad
\Xi_1 F'$.

From this equivalence we can get the needed biequivalence
$(G,G', \eta_3, \delta_3, \varepsilon_2, \zeta_2, \Phi_3, \Xi_3, \Psi_2,
\Omega_2):\T_2 \lra \T_3$
where

$G' \quad=\quad F\circ (G \circ F)'$

$\eta_3 \quad = \quad F \alpha \circ \zeta_1 : id_{\T_2} \lra FF'
\lra F\circ (G \circ F)' \circ G = G'G$

$\delta_3 \quad = \quad \varepsilon_1 \circ F \beta:G'G =  F\circ (G \circ F)' \circ G
\lra FF' \lra id_{\T_2}$

$\Phi_3 \quad = \quad F \Gamma \circ \Omega_1^{-1} : id
\lra \varepsilon_1 \circ \zeta_1 \lra \varepsilon_1 \circ F\beta \circ
F\alpha \circ \zeta_1 \lra \delta_3 \circ \eta_3$

$\Xi_3 \quad = \quad F \Theta \circ \Psi_1^{-1} : \eta_3 \delta_3
= F \alpha \circ \zeta_1 \circ \varepsilon_1 \circ F \beta \lra
F\alpha \circ F \beta \lra id.$ $\Box$

{\bf Retracts}. If $F$ is a retract of $G$, that is, if there exists a
commutative diagram as follows:

$$\xymatrix{ \T_1 \ar[rr]^H \ar[dd]_F \ar@/^1pc/[rrrr]^{id}&&
\T_3\ar[rr]^J\ar[dd]_G && \T_1\ar[dd]^F
\\
\\
\T_2 \ar[rr]_{H'} \ar@/_1pc/[rrrr]_{id}&& \T_4\ar[rr]_{J'}&&
\T_2}$$ and if $G$ is a weak equivalence (resp. fibration,
cofibration) then $F$ is also a weak equivalence (resp. fibration,
cofibration).

Weak equivalence. If $(G,G', \eta, \delta, \varepsilon, \zeta,
\Phi, \Xi, \Psi, \Omega):\T_3 \lra \T_4$ is a weak equivalence,
then so is $$(F,JG'H, J\eta H, J\delta H , J'\varepsilon H',
J'\zeta H' , J\Phi H , J\Xi H, J'\Psi H', J' \Omega H'):\T_1 \lra
\T_2.$$

Fibrations. Let $f \in_1 \T_1$ and $\beta: Ff \lrat g$ be a 2-cell in
$\T_2$. A lifting ``across'' $F$ of $f$ and $\beta$ will be denoted as
$L_F(f, \beta): f \lrat f'$. The needed lifting can then be described as
$$ L_F(f, \beta) \quad = \quad J(L_G(H(f), H'(\beta)))$$

Cofibrations. By assumption, $G, H$ and $H'$ are injective on 1-cells
and the left square commutes, therefore $F$ is injective on 1-cells.
$\Box$

{\bf First Lifting Axiom}. Consider the following commutative
diagram: $$\xymatrix{ \T_1\ar[rr]^U \ar[dd]_F && \T_3 \ar[dd]^G
\\
\\
\T_2 \ar[rr]_V \ar@{-->}[rruu]^H&& \T_4}$$ where $F$ if a
cofibration and $G$ is a fibration. The first axiom asserts that
if $F$ is also a weak equivalence, then there exists a lifting
$H$ making the two triangles commute.

Since $F$ is a trivial cofibration, $F$ is an inclusion of a full
sub-2-category. By Proposition 4(1), we can construct an $F':\T_2
\lra \T_1$ such that $F' \circ F = id_{\T_1}$ Define $H$ on
1-cells as follows: Let $f \in_1 \T_2$. Then $VFF'(f)$ is in the
image of $G$ and by the commutativity of the square is equal to
$GUF'(f)$. From the biequivalence of $F$ there is an iso-2-cell
$\gamma_f:GUF'(f) \lrat V(f)$. Since $G$ is a fibration, there is
a iso-2-cell in $\T_3$, $\delta_f:UF'(f) \lrat H(f)$. Use this as
a definition of $H$ on 1-cells.

Let $\alpha:f \lra f'$ be a 2-cell in $\T_2$. Then define $H$ on
2-cells as $$ H(\alpha)=\delta_{f'} \circ UF'(\alpha) \circ
\delta_f^{-1}: H(f) \lrat UF'(f) \lra UF'(f') \lrat H(f').$$ Such
an $H$ satisfies our requirements. $\Box$

{\bf Second Lifting Axiom}. Let $F$ be a cofibration and $G$ be a
fibration as in the previous axiom. The second lifting axiom
states that if $G$ is also a weak equivalence, then an $H$ exists
making the triangles commute.

By assumption, $F$ is injective on 1-cells and $G$ is surjective
on 1-cells. By a simple diagram chase, there is an $H:\T_2 \lra
\T_3$ defined on 0-cells and 1-cells. From the fact that G is a
biequivalence and Proposition 6, $G:\T_3(m,n) \lra \T_4(m,n)$ is
an equivalence. For any $f,f':m \lra n \in_1 \T_1$, we have from
the properties of an equivalence that $G:\T_3(H(f),H(f')) \lra
\T_4(V(f),V(f'))$ is an isomorphism. Use this isomorphism as a
definition for extending the lifting $H$ to 2-cells. $\Box$

{\bf First Factorization Axiom}. Every 2-theory-morphism $F:\T_1
\lra \T_2$ can be factored as a trivial cofibration $K:\T_1 \lra
\T_3$ followed by a fibration $G:\T_3 \lra \T_2$. $\T_3$ is the
categorical version of a path space.

The 0-cells of $\T_3$ are - of course - the natural numbers. It
suffices to describe the 1-cells of $\T_3$ as 0-cells of
$\T_3(1,m)$ for all $m\in \N$. The 0-cells of $\T_3(1,m)$ are
triples $(f,\alpha,g)$ where $f\in_0 \T_1(1,m)$, $g\in_0
\T_2(1,m)$ and $\alpha:Ff \lrat g \in_1 \T_2(1,m)$. 2-cells in
$\T_3$ are defined as $$ \T_3((f,\alpha, g), (f', \alpha', g'))
\quad = \quad \T_1(f,f').$$

$K:\T_1 \lra \T_3$ is defined as $$ K(f)=(f,id_{Ff},Ff) \qquad
K(\gamma:f \lra f')=\gamma:(f,id_{Ff},Ff) \lra
(f',id_{Ff'},Ff').$$ $G:\T_3 \lra \T_2$ is defined as $$
G((f,\alpha,g)) = g \qquad G(\gamma:(f,\alpha,g) \lra
(f',\alpha',g'))=\alpha' \circ G \gamma \circ \alpha^{-1}:g \lra
g' .$$

The factorization is clear. Obviously $K:\T_1(1,m) \lra
\T_3(1,m)$ is full, faithful and injective on 1-cells. For every
$(f,\alpha,g) \in_1 \T_3(1,m)$ we have the isomorphism $$
id_f:K(f)=(f,id_{Ff},Ff) \lrat (f,\alpha ,g).$$ So $K$ is dense
and hence a trivial cofibration. As for $G$ being a fibration, let
$(f,\alpha, g) \in_2 \T_3$ and $\gamma: g \lrat g' \in \T_2$. Then
$id_f:(f, \alpha,g ) \lrat (f, \gamma \circ \alpha, g')$ is a
lifting of $\gamma$.
 $\Box$

{\bf Second Factorization Axiom}. Every 2-theory-morphism $F:\T_1
\lra \T_2$ can be factored as a cofibration $K':\T_1 \lra \T_4$
followed by a trivial fibration $G':\T_4 \lra \T_2$. $\T_4$ is
the categorical version of a mapping cylinder.

The 0-cells of $\T_4$ are again - of course - the natural numbers.
It suffices to describe the 1-cells of $\T_4$ as 0-cells of
$\T_4(1,m)$ for all $m\in \N$. $$(\T_4(1,m))_0 \quad = \quad
(\T_1(1,m))_0 \coprod (\T_2(1,m))_0.$$ Warning: it is not
necessarily the case that $(\T_4(n,m))_0 \quad = \quad
(\T_1(n,m))_0 \coprod (\T_2(n,m))_0.$ The structure of
$\T_4(n,m)$ is generally more complex but can be calculated from
$\T_4(1,m)$ . The 1-cells of $\T_4(1,m)$ are given as

$$T_4(f,f')= \left\{ \begin{array} {l@{\quad:  \quad \mbox{ if
}}l} \T_2(F(f),F(f')) & f,f' \in_1 \T_1\\ \T_2(F(f),f') & f \in_1
\T_1, f' \in_1 \T_2\\ \T_2(f,F(f')) & f \in_1 \T_2, f' \in_1
\T_1\\ \T_2(f,f') & f,f' \in_1 \T_2\end{array} \right. $$

The cofibration $K':\T_1 \lra \T_4$ is described by $$K'(f)=f
\qquad K'(\alpha:f \lra f') = F(\alpha):F(f) \lra F(f')$$ The
fibration $G':\T_4 \lra \T_2$ is described by

$$G'(f)= \left\{ \begin{array} {l@{\quad:\quad \mbox{ if }}l}
 F(f) & f \in_1 \T_1\\
f & f \in_1 \T_2 \end{array} \right. $$ and $G(\alpha)=\alpha$.

The factorization is clear. The fact that $K'$ is a cofibration is
obvious. $G'$ is fibration because $G'$ has the (unique) lifting
property. $G$ is also surjective on 1-cells and locally full and
faithful. $\Box$

\begin{teo}
The category of 2-theories and 2-theory-morphisms admits a
functorial closed Quillen model category structure.
\end{teo}
By inverting the weak equivalences, we get the category
$Ho(\TTh)$ and the functor $\gamma: \TTh \lra Ho(\TTh)$ which
satisfies the universal property stated in the introduction..

\section{Universal Properties of Coherence}

The following proposition states the universal properties of the
mapping cylinder formed in the Second Factorization Axiom.
\begin{prop}
Let $F:\T_1 \lra \T_2$ and let $G \circ K : \T_1 \lra \T_4 \lra
\T_2$ be the factorization constructed in the second
factorization axiom. For all factorizations $G' \circ K': \T_1
\lra \T_5 \lra \T_2$ of $F$ consisting of a cofibration followed
by a trivial fibration, there is a unique {\bf isomorphism class}
of 2-theory-morphisms $H: \T_4 \lra \T_5$ making the triangles in
the following diagram commute:
\end{prop}

$$\xymatrix{ \T_1 \ar[rrrr]^F\ar[drr]^K\ar[rrddd]_{K'}&&&&\T_2 \\
&&\T_4\ar[urr]^G \ar@/_/[dd]_H^\Leftrightarrow \ar@/^/[dd]^{H'}
\\
\\&&\T_5. \ar[rruuu]_{G'}}$$

That is, for every factorization $G' \circ K' = F$ there is an $H
:\T_4 \lra \T_5$ such that $H \circ K = K'$ and $G' \circ  H =
G$. If there is any other $H'$ that satisfies these properties,
then for all $f \in_1 \T_4$ there is an iso-2-cell $\alpha_f :
H(f) \lrat H'(f)$ such that $\alpha \circ K = K'$ and $G' \circ
\alpha = G$. Note that all $H$ and $H'$s are trivial fibrations.

Remark: Proposition 8 is not a statement about all FCQMC
structures since, in general, model categories do not necessarily
have 2-cells.

{\bf Proof.} From the Second Lifting Axiom and from the fact that
the $G\circ K = F =G' \circ K'$ we have at least one $H$ making
the necessarily diagram commute. $$\xymatrix{ \T_1 \ar[rr]^{K'}
\ar[dd]_K && \T_5 \ar[dd]^{G'}\\
\\
\T_4 \ar[rr]_G \ar@{-->}[rruu]^H && \T_2.}$$

From the fact that $G' \circ H = G$, both $G$ and $G'$ are
biequivalences and the Two Out of Three Axiom, we have that $H$
is a biequivalence. Furthermore, let $G_1$ be the quasi-inverse of
$G'$. Then for all $f \in_1 \T_4$, $H(f)$ is isomorphic to $(G_1
\circ G)(f)$ say by $$\beta_f:H(f) \lrat (G_1\circ G)(f).$$

For any other $H':\T_4 \lra \T_5$ that satisfies the
commutativity of the triangles, we also have that $H'(f)$ is
isomorphic to $(G_1 \circ G)(f)$ say by $$\beta'_f:H'(f) \lrat
(G_1 \circ G)(f).$$ We can now define $$\alpha_f=(\beta'_f)^{-1}
\circ \beta_f : \quad H(f) \lrat (G_1 \circ G)(f) \lrat H'(f)$$
$\alpha_f$ satisfies the necessary requirements. $\Box$

Similar universal properties can be said for the First
Factorization Axiom.

In order to see Proposition 8 in action, let us work-out a
concrete example. let $\T_1=\Tbin$ be the 2-theory of anomic
multiplicative categories. That is, categories with a bifunctor
and no associating isomorphism assumed. $\Tbin(1,2)$ has one
element (the bifunctor) and $\Tbin(1,4)$ has five distinct
objects (the fourth Catalan number) with no morphisms between
them. Let $\T_2 = \Tsmon$ be the 2-theory of strict monoidal
categories. $\Tsmon(1,n)=\{*\}$ for all $n\in\N$ (we ignore
units). $\Tbin$ and $\Tsmon$ respectively represent the free-est
and strictest structures one can place on a category with a
bifunctor. Following the construction of the Second Factorization
Axiom, we get $\T_4 = \Tstar$ the 2-theory of starfish categories.

\begin{defi} A {\bf starfish category} is a category $\C$ with two
bifunctors $\otimes, \oplus : \C\times \C \lra \C$ such that

1) $A \oplus (B \oplus C) = (A \oplus B) \oplus C$

and

2) there exists a iso-natural transformation $\delta_{A,B}:
A\otimes B \lra A \oplus B$.
 \end{defi}

 For the association  $A\otimes (B \otimes C)$ we have the
 following naturality diagram:
 $$\xymatrix{ A\otimes (B \otimes C)
\ar[rr]^{\delta_{A, B \otimes C}} \ar[dd]_{id_A \otimes
\delta_{B, C}} &&A\oplus (B \otimes C) \ar[dd]^{id_A \oplus
\delta_{B, C}}\\
\\
A\otimes (B \oplus C) \ar[rr]^{\delta_{A, B \oplus C}} &&A\oplus
(B \oplus C).}$$ From the commutativity of this square, we have a
unique map of $\delta$s from $A\otimes (B \otimes C)$ to $A\oplus
(B \oplus C)$

By similar reasoning, we can extend this to all associations. For
any associated word $w$, let $w^\otimes$ (resp. $w^\oplus$)
represent the functors with only $\otimes$ (resp. $\oplus$)
between the letters. By naturality of $\delta$ there is a unique
$$\delta_w :w^\otimes \lra w^\oplus$$

$\Tstar(1,4)$ corresponds to the following: {\tiny $$\xymatrix{ &
A \otimes ((B \otimes C) \otimes D)\ar[rd]^\sim_\delta && (A
\otimes (B \otimes C)) \otimes D \ar[ld]^\sim_\delta \\ A \otimes
(B \otimes (C \otimes D)) \ar[rr]^\sim_\delta&& A \oplus B \oplus
C \oplus D && ((A \otimes B) \otimes C) \otimes
D\ar[ll]^{\sim}_{\delta} \\ &&(A \otimes B) \otimes (C \otimes
D). \ar[u]^\sim_\delta}$$ }

Hence the name ``starfish''. The definition of starfish-morphisms
and starfish-natural transformations are left to the reader. It
is not hard to see that every starfish category is
starfish-equivalent to a monoidal category and vice versa.

Let us return to Proposition 8. The 2-theory of starfish
categories is the 2-theory constructed in the Second Factorization
Axiom. The 2-theory of monoidal categories also satisfies the
Factorization Axiom. Putting all this together, we have
$$\xymatrix{ \Tbin \ar[rrrr]^F\ar[drr]^K\ar[rrddd]_{K'}&&&&\Tsmon
\\ &&\Tstar \ar[urr]^G \ar@/_/[dd]_H^\Leftrightarrow
\ar@/^/[dd]^{H'}
\\
\\&&\Tmon. \ar[rruuu]_{G'}}$$
The 1-cell in $\Tstar(1,4)$ corresponding to $A \oplus B \oplus C
\oplus D$ is the only ``loose end''. Every other 1-cell is forced
by the requirements of the commutativity of the triangles. The
multiplicity of choices where one can send that 1-cell correspond
to the multiplicity of the different $H$'s. However all the
different places where $H$ can take the ``loose ends'' are
(uniquely) isomorphic. Since for all $f \in_1 \Tstar$,
$\Tmon(H(f), H'(f))$ is isomorphic to $$\Tsmon(G'H(f), G'H'(f))
\quad = \quad \Tsmon(G(f),G(f)) \quad = \quad \{*\},$$ the
$\alpha$ of Proposition 8, is in fact unique.

$\Tstar$ is the free-est structure that can be added to $\Tbin$
and still be Morita equivalent to $\Tsmon$. $\Tmon$ is a type of
quasi-quotient of $\Tstar$ and hence also has this property. This
is a universal property of $\Tmon$.

There are other types of universal properties that can be said
about coherence from our point of view. We can place a FCQMC
structure on the category of 1-theories and 1-theory-morphisms.
The FCQMC structure is trivial and hence not very interesting in
itself. However, it interacts well with the FCQMC structure on
$\TTh$. For $\Th$, the weak equivalences are
1-theory-equivalences which are exactly 1-theory-isomorphisms.
The fibrations are all 1-theory-morphisms and the cofibrations
are all 1-theory morphisms that are injective on 1-cells. In both
model categories, all objects are fibrant and cofibrant.

We have the following Proposition from \cite{Quillen}.
\begin{prop}
Let $\C$ and $\C'$ be model categories and let $$\xymatrix{ \C
\ar@/^1pc/[rr]^L_\bot && \C' \ar@/^1pc/[ll]^R}$$ be a pair of
adjoint functors. Suppose $L$ preserves cofibrations and $R$
preserves fibrations. Then the left Kan extension
$Lan_\gamma(\gamma' \circ L)$ and the right Kan extension
$Ran_{\gamma'}(\gamma \circ R)$ exists and are adjoint.
\end{prop}

 $$\xymatrix{ \C \ar@/^1pc/[rr]^L_\bot \ar[dd]_\gamma
&& \C' \ar@/^1pc/[ll]^R \ar[dd]_{\gamma'}\\ \\ Ho(\C)
\ar@/^1pc/[rr]^{Lan_\gamma(\gamma' \circ L)}_\bot && Ho(\C')
\ar@/^1pc/[ll]^{Ran_{\gamma'}(\gamma \circ R) } }$$

Placing this in our context, only the $U \vdash d$ adjunction of
Diagram (\ref{XY:adj}) satisfies the requirements of Proposition
9. And so we have $$\xymatrix{ \Th \ar@/^1pc/[rr]^d_\bot
\ar@{=}[dd]_{id} && \TTh \ar@/^1pc/[ll]^U \ar[dd]_{\gamma'}\\ \\
Ho(\Th) \ar@/^1pc/[rr]^{\gamma' \circ d}_\bot && Ho(\TTh).
\ar@/^1pc/[ll]^{Ran_{\gamma'}(U)}  }$$

The other adjunctions between $\Th$ and $\TTh$ do not satisfy the
requirements of the Proposition 9. Nor do any of them induce an
equivalence between $\Th$ and $Ho(\TTh)$. This does not prove
that no such equivalence exists, but we believe the two categories
are, in fact, not equivalent (How does one prove two categories
are {\it not} equivalent?)

We would like to point to places in the literature that seem to
be instances of the $\gamma' \circ d$ functor. Fangjun Arroyo
\cite{Arroyo} has proven that symmetric monoidal categories are
precisely the homotopy commutative monoids in $\Cat$ where the
weak equivalence in $\Cat$ are equivalences of categories. In our
language this, in effect, becomes $\gamma' \circ d(T_{com-mon}) =
[\Tsym]$ where $T_{com-mon}$ is the 1-theory of commutative
monoids and $[\Tsym]$ is the homotopy class of the 2-theory of
symmetric monoidal categories.

Recently, Tom Leinster \cite{Leinster} has proven a similar result
for monoids and monoidal categories. We are left with the obvious
question of where does braided monoidal categories fit in this
scheme?

We would like to conclude by stating that if one assumes that the
Kronecker bifunctor defined in Section 2 is symmetric ($\T_1 \KP
\T_2 \cong \T_2 \KP \T_1)$, then it extends to the homotopy
category. All we have to do is show that $\KP$ takes two
biequivalences to a biequivalence. This is a short lemma if one
takes into account Proposition 3 and the way that $F_1 \KP F_2$
is defined. And hence we have the following diagram which will
help us build new coherence theorems from old ones. $$\xymatrix{
\TTh \times \TTh \ar[rr]^\KP \ar[dd]_{\gamma \times \gamma} &&
\TTh
\ar[dd]^\gamma \\
\\
Ho(\TTh)\times Ho(\TTh) \ar[rr]^{Ho(\KP)} && Ho(\TTh). }$$

\begin{examp} In \cite{Paper3} we have shown that $$\Tsmon \KP \Tsmon
\cong \Tsbraid .$$ In Section 3, we have shown that $\Tsmon$ is
Morita equivalent to $\Tmon$ and $\Tsbraid$ is Morita equivalent
to $\Tbraid$. From the above commutative square, we see that
$$\Tmon \KP \Tmon \cong \Tbraid.$$
\end{examp}

\section{Future Directions}

{\bf The $\Th$, $Ho(\TTh)$ Adjunction.} We have shown that the
usual adjunctions between $\Th$ and $\TTh$ do not induce an
equivalence between $Ho(\Th)=\Th$ and $Ho(\TTh)$. However, this
does not mean that no such equivalence exists. Although we
conjecture that the categories are in fact not equivalent, we
have, as yet, no idea of how to prove this. It is left as an open
question. One must realize that this question goes against the
entire grain of the subject. Quillen invented model categories to
show when two different model categories have the same homotopy
category. We are asking to show that the different model
categories have different homotopy categories.

A moment of speculative thought is in order. Assume that $\Th$ is
not equivalent to $Ho(\TTh)$. This would show that coherence is
not simply a homotopical notion. Rather it is also in an
intrinsic manner an algebraic notion. The adjunction discussed
after Proposition 9 is - as all adjunctions - an algebraic
concept. There are generators, relations, free and forgetful
functors, universal properties etc. Coherence seems to be a
complex notion which encompasses both elements of homotopy {\bf
and} elements of algebra.

{\bf The Structure -- Semantics Adjunction.} Between any two
$n$-categories there are $n+1$ different notions of equivalence
that can connect them \cite{Street}. At the present time, we are
looking at the diverse homotopy categories these different
equivalences induce on the category of ${\bf nCat}$ \cite{Paper4}.
In particular, we shall look at the category of $\TCat$. More to
the point, we plan on examining the subcategory of $\TCatcat$
where the semantics of algebraic structure lives \cite{Paper3,
Lawvere}. Our goal will be to determine the extent to which the
Structure-Semantics (quasi-)adjunction preserves the Quillen
model category structures of $\TTh$ and $\TCatcat$.

{\bf Higher Cells of $\TTh$.} By only looking at the (1-)category
of 2-theories and 2-theory-morphisms, we are ignoring the higher
cells of $\TTh$. $\TTh$ is a 3-category, but we have effectively
disregarded 2-theory-natural transformations and
2-theory-modifications (although they are used surreptitiously in
the characterization of 2-theory-biequivalences.) Abandoning the
higher structure seems unnatural. There are many other model
categories having higher structure that should not be ignored.
Surely there are important theorems that can be proved about these
higher cells. To our knowledge, no one has written down axioms
for a Quillen model 2-category or 3-category. How should
(co)fibrations behave with respect to (iso-)2-cells? How is the
fraction category constructed when there are higher cells
involved? etc. The task of writing down such an axiom system is
far beyond the author's capabilities. We are simply pointing to a
glaring gap in the literature with the hope that someone takes on
the challenge .

{\bf Relative Homotopy Theory.} One of the central themes in
coherence theory is that when dealing with morphisms between
algebras, certain operations should be preserved up to an
(iso)morphism and certain operations are preserved strictly.
\cite{Paper3} dealt with this by having a controlling 2-theory
$\T_1$ and a 2-theory-morphism $G:\T_1 \lra \T_2$ that decides
what type of preservation property an operation in $\T_2$ should
have. If $\T_1$ controls the operations in $\T_2$ and $\T_3$,
then a biequivalence between $\T_2$ and $\T_3$ should be strict
on $\T_1$
$$
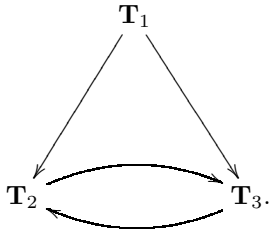
\xymatrix{& \T_1 \ar[ldd] \ar[rdd] \\
\\
\T_2 \ar@/^1pc/[rr] && \T_3. \ar@/^1pc/[ll]}$$ This seems very
similar to relative homotopy theory where one has a subspace
(subdiagram) such that a homotopy is the identity on the
subspace. Alex Heller \cite{HellerRel} has given axioms (in the
same language as \cite{HellerTheory}) for relative homotopy
theory. We are planning to make the connection between relative
homotopy theory and our study of coherence.

{\bf Generalizations of 2-Theories.} 2-Theories can not describe
all structures one usually places on a category. In particular,
we must extend the definition of a 2-theory to handle
contravariant functors and dinatural transformations. A method of
doing this was discussed in \cite{Paper3}. If we are successful
in formulating the notion of a generalized 2-theory, then we will
be able to describe closed categories and hence enter the world
of low-dimensional topology, quantum groups and computer science
(e.g. traced monoidal categories). It would be nice to extend the
results in this paper to generalized 2-theories. Perhaps we will
be able to place a FCQMC structure on the category of 2-monads
\cite{Blackwell}. Upon entering the world of, say, low
dimensional topology, we might ask what it means for one
topological invariance to be of the same homotopy type as
another? Similar questions in other areas are very interesting.

Bloom {\it et al} \cite{Bloom} has extended 2-theories in another
direction. They have defined iteration 2-theories These are
2-theories with extra operations that are useful in describing
feedback and fixed points. See \cite{Wagner} for a survey of many
such interesting 1-theories and 2-theories. Such generalizations
is of extreme importance to computer science. They are used in
describing rewrite systems, trees, data types, etc. Our goal is
to extend this work to incorporate iteration 2-theories. We hope
to answer questions as to when two data types are ``of the same
homotopy type''? When do two rewrite systems produce the same
language ``up to homotopy''?

{\bf Algebraic Operads.} Vladimir Hinich \cite{Hinich} has placed
a closed Quillen model category structure on the category of
differential graded operads over a ring. Such operads are ways of
describing algebraic structures on chain complexes of modules. One
of the main ideas of quantum groups is the structure of an
algebra (coalgebra, bialgebra, Hopf algebra, quasi-Hopf
quasi-triangular algebra etc) A is reflected in the structure of
the category of modules (comodules, bimodules, bicrossed modules
etc) of A (see e.g. \cite{Kassel} or \cite{Yetter}.) Hence there
is some type of functor $Rep$ from the category of differential
graded operads to $\TTh$  that takes an operad $\O$ to the
2-theory of the structure of the category of modules from an
arbitrary $A\in Alg(\O)$. Let us explain. There are three levels
of algebraic structure here. There is (i) a category of operads,
$OPERADS$; (ii) for each operad $\O \in OPERADS$, there is a
category of algebras/models of A, $Alg(\O)$; and (iii) for each
$A\in Alg(\O)$ there is a category of modules of $A$, $Mod(A)$. An
operad in $OPERAD$ determines the type of structure of (iii).
Types of structures in (iii) are described by 2-theories. So we
have a functor from $OPERADS$ to $\TTh$.

Questions: Can we formally describe this functor $Rep$? Is there
an inverse (quasi-adjoint, adjoint) of $Rep$? Does $Rep$ respect
Hinich's model structure? Will an inverse respect our model
category structure? What is the relationship of homotopy theory to
representation theory?

{\bf Other Notions of Weak Equivalence.} Not every coherence
relationship of interest is a biequivalence. $\Tmon$ and $\Tsmon$
are biequivalent. However, the relationship between $\Tbraid$ and
$\Tsym$ is not so simple. We believe that they are
quasi-biequivalent. A quasi-biequivalence is a weakening of the
concept of a biequivalence where we substitute a
2-theory-quasi-natural transformation instead of a
2-theory-natural transformation. We would like to investigate
other FCQMC structures on $\TTh$ then the one given in this
paper. Different sets of weak equivalences give different
homotopy categories and hence different notions of coherence.

{\bf Tools of Homotopy Theory} Once there is a FCQMC structure on
a category $\C$ one can explore and exploit the structure of $\C$
with the powerful tools of homotopy theory. We might look at
homotopy limits and colimits, homotopy Kan extensions, long exact
sequences; homology etc. We plan on going on and looking at
$\TTh$ with these tools. Some further questions arise: Are there
``minimal models'' of a homotopy class of 2-theories? Although
there are no Postnikov towers for $\TTh$ (it is not a {\it
pointed} FCQMC), can we nevertheless decompose algebraic
2-theories ``up to homotopy''. Much work remains to be done.

Department of Computer and Information Science\\
Brooklyn College, CUNY\\
Brooklyn, N.Y. 11210\\
email: noson@sci.brooklyn.cuny.edu
\end{document}